\documentclass[12pt]{article}

\usepackage{amssymb}
\usepackage{amsmath}
\usepackage{mathdots}
\usepackage{amsbsy}
\usepackage{amscd}
\usepackage{amsthm}
\usepackage{mathrsfs}
\usepackage{verbatim}
\usepackage[colorlinks]{hyperref}
\usepackage{fullpage}

\usepackage{url}

\usepackage[english]{babel}
\usepackage[T1]{fontenc}
\usepackage[utf8]{inputenc}
\usepackage{tikz}
\usepackage{authblk}

\usepackage{multirow,tabularx}

\theoremstyle{definition}
\newtheorem{Theorem}{Theorem}[section]
\newtheorem{Lemma}[Theorem]{Lemma}

\newtheorem{Remark}[Theorem]{Remark}

\newtheorem{Notation}[Theorem]{Notation}

\newcommand{\C}{\mathbb{C}}

\newcommand{\Q}{\mathbb{Q}}

\newcommand{\mt}[1]{\text{#1}}

\begin{document}

\title{Diagonal Orbits in a Type A Double Flag Variety 
of Complexity One}
\author[1,*]{Mahir Bilen Can}
\author[2]{Tien Le}
\affil[1]{{\small Tulane University, New Orleans, USA}}
\affil[2]{{\small Tulane University, New Orleans, USA}}
\affil[*]{{\small Corresponding author: Mahir Bilen Can, mahirbilencan@gmail.com}}
\normalsize

\date{\today}
\maketitle

\begin{abstract}
We continue our study of the inclusion posets of diagonal $SL(n)$-orbit 
closures in a product of two partial flag varieties.
We prove that, if the diagonal action is of complexity one,
then the poset is isomorphic to one of the 28 posets that 
we determine explicitly. Furthermore, our computations show
that the number of diagonal $SL(n)$-orbits in any
of these posets is at most 10 for any positive integer $n$.
This is in contrast with the complexity 0 case, where, in some cases, the resulting posets 
attain arbitrary heights. 
\vspace{.5cm}

\noindent 
\textbf{Keywords:} Double flag varieties, 
complexity one actions, lattices.\\ 
\noindent 
\textbf{MSC: 06A07, 14M15} 
\end{abstract}

\section{Introduction}\label{S:Introduction}

Let $G$ be a connected reductive complex algebraic group, 
and let $B$ be a Borel subgroup in $G$. 
Let $X$ be an irreducible complex algebraic $G$-variety. 
We denote the action of $G$ on $X$ by $G:X$. 
A typical example for such a variety is the 
homogeneous space $G/H$, where $H$ is a closed subgroup of $G$,
and the action of $G$ on $G/H$ is given by the 
multiplication action of $G$ on the left cosets of $H$ in $G$.
The complexity of $G:X$, denoted by $c_G(X)$, 
is defined as the codimension of a general $B$-orbit in $X$.
This notion plays an important role in the theory of equivariant  
embeddings of homogeneous spaces, see~\cite{Timashev}.
As it is demonstrated by the seminal paper of Panyushev~\cite{Panyushev99}, 
among all homogeneous spaces of $G$, the ones with complexity at most one form the most remarkable subclass.

An enduring problem in representation theory is 
to decompose the tensor products of irreducible 
representations of $G$. 
Let $\lambda_i$ ($1\leq i \leq 2$) be two dominant 
weights corresponding to the irreducible representations 
$V_i$ ($1\leq i \leq 2$) of $G$, 
and let $P_i$ ($1\leq i \leq 2$) denote the 
corresponding parabolic subgroups 
that arise as the stabilizer subgroups of highest weight 
vectors $v_i\in V_i$ ($1\leq i \leq 2$).
There is a close relationship between the decomposition of 
$V_1 \otimes V_2$ as
a $G$-module and the polynomial invariants of the  
diagonal action of $G$ on the double flag variety 
$X:=G/P_1\times G/P_2$.
By using the coordinate ring of the 
affine cone over the double flag variety, 
in~\cite{Littelmann}, Littelmann obtained 
precise description of the 
decompositions of the tensor products of two 
fundamental representations of simple groups.
This progress motivated the works~\cite{MWZ1,MWZ2}, 
and~\cite{Stembridge}.
In the last reference, Stembridge classified all 
multiplicity-free tensor products 
of irreducible representations of semisimple complex 
Lie groups. This classification amounts to the classification of the parabolic subgroups 
$P_i$ ($1\leq i\leq 2$) such that $c_G(G/P_1\times G/P_2)=0$. 
Finally, in~\cite{Ponomareva}, Ponomareva classified all 
double flag varieties of complexity one.
In the same paper, Ponomareva showed by examples 
how one could use the results of 
Brion~\cite{Brion89} and Timashev~\cite{Timashev00} 
for decomposing the 
spaces of global sections of the line bundles on a double 
flag variety of complexity $\leq 1$.
In the present article, we focus on the double flag varieties of complexity one.
Our purpose here is to give a complete description of the inclusion order on the closures of the $G$-orbits 
in $G/P_1\times G/P_2$ when $c_G(G/P_1\times G/P_2)=1$.

To further motivate our discussion, 
let us mention another setup where the diagonal orbits are of crucial importance. 
In~\cite{DeligneLusztig}, Deligne and Lusztig constructed
the complex linear representations of finite groups of Lie type by using the $\ell$-adic cohomology with compact support
on certain varieties. Let $\mathbf{F}_q$ denote the finite field with $q$ elements,
let $G$ be a reductive group defined over an algebraic closure of $\mathbf{F}_q$,
and let $F$ denote a Frobenius map on $G$. 
Let $w$ be an element of the Weyl group of $G$. 
The {\em Deligne-Lusztig variety associated with $w$}, denoted by $X(w)$, 
consists of all Borel subgroups $B$ of $G$ such that $B$ and $F(B)$ are in relative position $w$. 
In other words, $X(w)$ is the intersection of the $G$-orbit corresponding to $w$ in 
$G/B\times G/B^-$, where $B^-$ is the unique opposite Borel subgroup to $B$, with the graph of the Frobenius map.
More recently, Digne and Michel extended this theory to the setting of partial flag varieties, see~\cite{DigneMichel}.
In essence, the poset that we study in our paper is about the natural hierarchy between 
the parabolic Deligne-Lusztig characters, namely, the characters of the representations of $G^F$ on 
$H^*_c(X(w),\overline{\Q_\ell})$, where $X(w)$ is a Deligne-Lusztig variety in $G/P\times G/P^-$.
Here, $P^-$ is a parabolic subgroup such that $P\cap P^-$ is a common Levi
subgroup of both of $P$ and $P^-$.
Finally, let us mention that the same partial order arises rather naturally in the study of the nilpotent 
variety of the dual canonical monoids, see~\cite{Therkelsen_Thesis} and~\cite{Can19}.

Let $X$ be a normal $G$-variety, and let $B$ denote a Borel subgroup of $G$. 
In many ways the geometry of $X$, as a $G$-variety, 
depends on how $G$- and $B$-orbits in $X$ fit together. 
For example, if $X$ has finitely many $G$-orbits, 
 then the rational Chow group of $X$ has a decomposition with respect to $G$-orbits, see~\cite{BinghamCanOzan}.
With this fact in mind, in our earlier work~\cite{Can18}, for $G=SL(n)$, 
we showed that if $c_G(X)=0$, 
then the inclusion poset of $G$-orbit closures in $X$ 
is a particular kind of graded lattice; it is either a chain,
or it is what we called a `ladder poset.' 
In higher complexity, these posets 
can be very complicated; they are not necessarily graded. 
However, they always have a unique minimal and a unique 
maximal element.
In the case of complexity one, as we show, 
most of them turn out to be lattices, and not all of them are graded.
Our main theorem is the following statement.

\begin{Theorem}\label{T:prelim}
Let $G$ denote $SL(n)$ and let $X$ be a double 
flag variety $G/P_1\times G/P_2$. If $c_G(X)=1$,
then the inclusion poset of $G$-orbit closures 
in $X$ is one of the 28 posets whose Hasse diagrams 
are as in Figure~\ref{F:28}.
\end{Theorem}

For us, the most surprising outcome of our computation 
is the number of $G$-orbits in $G/P_1\times G/P_2$.
Although there are infinitely many complexity one 
double flag varieties, in each case, the number of $G$-orbits turns out
to be bounded by 10; this is in contrast with the complexity zero case, where 
there are infinitely many non-isomorphic $G$-orbit containment posets, 
and they can be of arbitrary height.

Next, we give a brief description of our paper. 
In Section~\ref{S:Preliminaries},
we present some background material regarding 
our posets. Section~\ref{S:Complexity} forms 
the main body of our paper; we depict the Hasse 
diagrams of our posets 
in Figure~\ref{F:28}.
The subsequent Section~\ref{S:Conclusions} 
is the concluding section for the proof of our main theorem.  
Finally, in Section~\ref{S:Final}, we mention an alternative method
for proving our theorem.

\section{Preliminaries}\label{S:Preliminaries}

\subsection{}\label{SS:1}

Let $G$ be a complex semisimple algebraic group,
let $B$ be a Borel subgroup in $G$, and let $T$ 
be a maximal torus of $G$ that is contained in $B$. 
We denote by $\Phi$ the root system corresponding
to the pair $(G,T)$, and we denote by $\Delta$ the set of 
simple roots determined by $B$. 
A parabolic subgroup $P$ of $G$ is said to be 
{\em standard with respect to $B$} if $B\subseteq P$.  
In this case, $P$ is uniquely determined by a subset 
$I \subseteq \Delta$ such that $|I| = \dim P/B$. 

Let $N_G(T)$ denote the normalizer subgroup of 
$T$ in $G$. 
The {\em Weyl group} $W:=N_G(T)/T$ of $G$
is a Coxeter group,
and we denote its Coxeter generating system 
corresponding to $\Delta$ by 
\[
R(\Delta) :=\{ s_\alpha \in W :\ \alpha \in \Delta \}.
\]
The elements of $R(\Delta)$ are called 
the {\em simple reflections} relative to $B$.
If the set of simple roots we are using 
is fixed, then we will denote $R(\Delta)$ 
by $R$ to ease our notation.

We will interchangeably use the letters $I$ and $J$ to denote  
subsets of $\Delta$ and the corresponding
subsets of simple reflections in $R(\Delta)$. 
The {\em length} of an element $w\in W$,
denoted by $\ell(w)$, is the 
minimal number of simple reflections 
$s_{\alpha_i}\in R(\Delta)$ that is needed for the equality  
$w=s_{\alpha_1}\cdots s_{\alpha_k}$ hold true.
In this case, the product $s_{\alpha_1}\cdots s_{\alpha_k}$
is called a {\em reduced expression} for $w$.

The {\em Bruhat-Chevalley order} on $W$ 
is be defined by declaring 
$v\leq w$ $(w,v\in W)$ 
if a reduced expression of $v$ 
is obtained from a reduced 
expression $s_{\alpha_1}\cdots s_{\alpha_k}=w$ 
by deleting some of the simple reflections $s_{\alpha_i}$
in $w$. More geometrically, the Bruhat-Chevalley order
is given by 
$v \leq w \iff B\dot{v} B/B \subseteq \overline{B \dot{w} B/B}$.
Here, $\dot{v}$ and $\dot{w}$ are any representatives
of $v$ and $w$ in $N_G(T)$, respectively. The 
sets $B\dot{v}B/B,B\dot{w}B/B$ denote the 
$B$-orbits of $\dot{v},\dot{w}$ in $G/B$, and 
the bar on $B\dot{w}B/B$ indicates the Zariski closure. 
In this notation, $\ell(w)$ is equal to the dimension
of the orbit $B\dot{w}B/B$.

Let $G$ be a classical semisimple matrix group with entries in $\C$,
and let $B$ denote its Borel subgroup consisting of upper triangular matrices.
The parabolic subgroups of $G$ containing $B$ 
have block-triangular structure, and they are determined by the 
sizes of the diagonal blocks. Following Ponomareva's notation
from~\cite{Ponomareva}, if $P$ is a parabolic subgroup 
containing $B$, then we will denote by 
$Bl(P)$ the sequence $(p_1,\dots, p_r)$, where $p_i$ 
denotes the size of the $i$-th block in $P_I$. 
For example, if $P$ is the Borel subgroup 
of upper triangular matrices in $SL(n)$, then 
each diagonal block of $P$ is a $1\times 1$ matrix,
therefore, $Bl(P_I)$ is the sequence $(1,1,\dots, 1)$ 
with $n$ entries.

Our primary example is the matrix group $G=SL(n)$. 
We take $B$ as the Borel subgroup of upper 
triangular matrices, and we take $T$ as the 
maximal torus of diagonal matrices in $B$. 
The Weyl group $W$ of $SL(n)$ is denoted by $S_n$,
which is isomorphic to the symmetric group 
of permutations of $\{1,\dots, n\}$. 
The set of simple roots relative to $B$, 
that is $\Delta_{n-1} := \{\alpha_1,\dots, \alpha_{n-1} \}$, 
is ordered so that the $i$-th simple reflection 
$s_{\alpha_i}$ ($1\leq i \leq n-1$) is the 
simple transposition $s_i\in S_n$ that interchanges 
$i$ and $i+1$. Thus we set 
$$
R_{n-1} := R(\Delta_{n-1}) = \{s_1,\dots, s_{n-1}\}.
$$
If a permutation
$w$ in $S_n$ is given in one-line notation 
$w=w_1\dots w_{n}$, then its length is equal to 
the cardinality of the following set:
$\{ 1\leq i < j \leq n :\ w_i > w_j \}$.

An important fact that we repeatedly use in our paper is that
$SL(n)$ is the stabilizer subgroup in $SL(n+1)$ 
of the standard basis vector $e_{n+1}$ of  
$\C^{n+1}$, where $SL(n+1)$ acts by its defining representation.
In particular, by using this identification of $SL(n)$ 
as a subgroup of $SL(n+1)$,  
we will use the following containments 
without further mentioning in the sequel: 
$$
\Delta_{n-1} \hookrightarrow \Delta_n,\ 
R_{n-1} \hookrightarrow R_n,\ 
\text{ and $S_n \hookrightarrow S_{n+1}$ (as a subgroup)}.
$$

\subsection{}\label{SS:2}

Let $X_1$ and $X_2$ be two $G$-varieties. 
Let $x_i \in X_i$ ($1\leq i \leq 2$) be two points 
in general positions. If $G_i \subset G$ denotes 
the stabilizer subgroup of $x_i$ in $G$, then 
$\mt{Stab}_G(x_1\times x_2)$ coincides with 
the stabilizer in $G_1$ of a point in general 
position from $G/G_2$ (or, equivalently, with 
the stabilizer in $G_2$ of a point in general position 
from $G/G_1$), see~\cite{Panyushev99}. 
As a special case, we consider 
the $G$-variety $X:= G/P_1\times G/P_2$. 
The proof of the following 
lemma is not difficult, see~\cite[Lemma 2.1]{Can18}.

\begin{Lemma}\label{L:poset isom}
The poset of $G$-orbit closures in $X$ is 
isomorphic to the poset of $P_2$-orbit closures
in $G/P_1$.
\end{Lemma}

From now on we assume that $P_1$ and $P_2$ are 
standard parabolic subgroups with respect to $B$. 
If $I$ and $J$ are the subsets of $R:=R(\Delta)$ 
(or, of $\Delta$) 
that determine $P_1$ and $P_2$, respectively, 
then we will write $P_I$ (resp. $P_J$) in place
of $P_1$ (resp. $P_2$).
The Weyl groups of $P_I$ and $P_J$ 
are denoted by $W_I$ and $W_J$, respectively. 
The set of 
$(W_I,W_J)$-double cosets in $W$ is 
denoted by $W_I \backslash W / W_J$.

\subsection{}

It follows from Bruhat-Chevalley decomposition that 
the set of $B$-orbits in $G/P_J$ are in 
a bijection with the set of minimal length left 
coset representatives 
for $W/W_J$, which we denote by $W^J$. 
The set of minimal length right 
coset representatives for $W_I\backslash W$ is denoted
by $^{I}W$. 
In a similar way, $W_I \backslash W / W_J$
is in a bijection with the set of $P_I$-orbits in $G/P_J$,
see~\cite[Section 21.16]{Borel}.
Let $w$ be an element from $W$,
and let $[w]$ denote the double coset $W_I w W_J$. 
Let 
$$
\pi : W \rightarrow W_I \backslash W / W_J
$$ 
denote the canonical projection onto the set of 
$(W_I,W_J)$-double cosets. 
Then the preimage in $W$ of every double coset 
in $W_I\backslash W / W_J$ is an interval with 
respect to Bruhat-Chevalley order.
Therefore, there is a unique maximal and a 
unique minimal element, see~\cite{Curtis85}.
Moreover, if $[w]$ and $[w']$ are two elements from $W_I\backslash W / W_J$, 
and $\bar{w}$ and $\bar{w'}$ are the maximal elements in the cosets $[w]$ and $[w']$, respectively, 
then $w \leq w'$ if and only if $\bar{w} \leq \bar{w'}$.
(This can be seen directly by a geometric argument, but see~\cite{HohlwegSkandera} also.)
Therefore, the set of $(W_I,W_J)$-cosets has a natural combinatorial partial order 
defined by 
$$
[w] \leq [w'] \iff w \leq w' \iff \bar{w} \leq \bar{w'}, 
$$
where $[w],[w'] \in W_I\backslash W / W_J$.
There is a geometric interpretation of this partial order:
If $O_1$ and $O_2$ are two $P_I$-orbits in $G/P_J$
with the corresponding double cosets $[w]$ and $[w']$,
respectively, then $O_1 \subseteq \overline{O_2}$ if and only if $w \leq w'$. 
The bar on $O_2$ stands for the Zariski closure in $G/P_J$.

Let $[w]$ ($w\in W$) be an element from 
$W_I\backslash W /W_J$ such that 
$\ell(w) \leq \ell(v)$ for all $v\in [w]$. 
Such minimal length double coset 
representatives are parametrized by the set 
${}^{I}{W} \cap W^J$. 
From now on, we denote $^{I}{W} \cap W^J$ by 
$U_{I,J}^-$. Set $H= I \cap w J w^{-1}$. 
Then $ uw \in W^J$ for $u\in W_I$ 
if and only if $u$ is a minimal length coset 
representative for $W_I/W_H$. 
In particular, every element of $W_I w W_J$ 
has a unique expression of the form 
$uwv$ with $u\in W_I$ is a minimal length 
coset representative of $W_I/W_H$, $v\in W_J$ and 
\begin{align}\label{A:additive length}
\ell(uwv) = \ell(u)+\ell(w)+\ell(v).
\end{align}

For $i\in \{1,\dots, n-1\}$, let $s_i$ denote the $i$-th simple transposition. 
Let $w$ be a permutation in $S_n$, and let $w_1\dots w_n$ be the one-line notation for $w$.
The number $i$ is called a {\em right descent of $w$} if $w_i > w_{i+1}$. 
Equivalently, $i$ is a right descent if $\ell(w s_i) < \ell(w)$. 
The set of all right descents of $w$, denoted by $\mt{Des}_R(w)$, is called the {\em right descent set of $w$}.
In a similar way, the integer $i$ is said to be a {\em right ascent of $w$} if $w_i < w_{i+1}$, 
or, equivalently, $\ell(ws_i ) > w$. 
The {\em right ascent set of $w$}, denoted by $\mt{Des}_R(w)$, is the set of all right ascents of $w$.
In this notation, the following characterization of  $U_{I,J}^-$ will be useful for our purposes: 
\begin{align*}
U_{I,J}^- &= \{ w\in W:\ I\subseteq \mt{Asc}_L(w)\ 
\text{ and } J\subseteq \mt{Asc}_R(w) \}\\
&=  \{ w\in W:\ I^c\supseteq \mt{Des}_R(w^{-1})\ 
\text{ and } J^c\supseteq \mt{Des}_R(w) \}
\end{align*}

\begin{Remark}\label{R:involution}
Let $\theta$ denote the involution of the set $R_{n-1}$ 
that is defined by $s_i \mapsto s_{n-i}$ for $i\in \{1,\dots, n-1\}$. 
Then $U_{I,J}^-$ and $U_{\theta(I),\theta(J)}^-$ are isomorphic as posets. 
\end{Remark}

\section{Computations}\label{S:Complexity}

As we mentioned before, Ponomareva~\cite{Ponomareva} has determined
the parabolic subgroups $P_I$ and $P_J$ in a semisimple complex algebraic
group $G$ such that the complexity of the diagonal action of $G$ on $G/P_I\times G/P_J$
is one. For $G=SL(n)$, the possible $P_I$ and $P_J$'s,
according to their block sizes, are listed in Table~\ref{T:Ponomareva}.
There are in total eight major cases.

\begin{table}[htp]
\begin{center}
  \begin{tabular}{| c | c | c | c | }
    \hline
 &    Number of blocks &    $Bl(P_I)$ & $Bl(P_J)$ \\ \hline
 1.&   $2,3$  & $(3,p_2),\ p_2 \geq 3$ & $(q_1,q_2,q_3),\ q_1, q_2, q_3 \geq 2$\\   \hline
 2. &  $2,3$ & $(p_1,p_2),\ p_1,p_2 \geq 3$ & $(2,2,q_3),\ q_3 \geq 2$\\   \hline
 3. & $2,3$  & $(p_1,p_2),p_1,p_2 \geq 3$ & $(2,q_2,2),\ q_2 \geq 2$\\   \hline
4. & $2,4$ & $(2,p_2),\ p_2\geq 3$ & $(q_1,q_2,q_3,q_4)$\\  \hline
5. & $2,4$ & $(p_1,p_2),\ p_1,p_2 \geq 2$ & $(1,1,1,q_4)$\\   \hline
6. & $2,4$ & $(p_1,p_2),\ p_1,p_2 \geq 2$ & $(1,1,q_3,1),\ q_3\geq 2$\\  \hline
7. &    $3,3$     &$(1,p_2,1)$, $p_2 \geq 2$ & $(q_1,q_2,q_3)$\\ \hline
8. &  $3,3$      &$(1,1,p_3)$, $p_3 \geq 2$ & $(q_1,q_2,q_3)$\\  \hline
  \end{tabular}
\end{center}
\caption{The list of all complexity $1$ double flag varieties for $G=SL(n)$.}
\label{T:Ponomareva}
\end{table}

{ 
In the rest of this section, we will describe the structure of the 
poset $U_{I,J}^-$ for each pair of parabolic subgroups $(P_I,P_J)$ 
from Ponomareva's list.
For $i\in \{1,\dots, 8\}$, the $i$-th row of Table~\ref{T:Ponomareva} will be analyzed in Subsection 3.i.}

\begin{Notation}
If $n$ is a positive integer, then we will use the shorthand $[n]$ to 
denote the set $\{1,\dots, n\}$. 
\end{Notation}

\subsection{$Bl(P_I)=(3,p_2)$, $p_2\geq 3$ and $Bl(P_J)=(q_1,q_2,q_3)$,
$q_1,q_2,q_3\geq 2$.} 

Let $n$ denote $3+p_2$, which is equal to $q_1+q_2+q_3$.
Clearly, $n \geq 6$ and $p_2 > q_3$. 
Since $I^c=\{s_3\}$, and $J^c=\{s_{q_1},s_{q_1+q_2}\}$,
we see that if $w=w_1\dots w_n \in U_{I,J}^-$, then 
\begin{enumerate}
\item[(i)] for $i\in [n-1]\setminus \{3\}$, $i$ comes before $i+1$ in $w$;
\item[(ii)] $w_1<\dots < w_{q_1}$, $w_{q_1+1}<\dots < w_{q_1+q_2}$,
$w_{q_1+q_2+1}<\dots < w_{n}$.
\end{enumerate}
This implies that $1\in \{w_1,w_{q_1+1},w_{q_1+q_2+1}\}$,
and that $n\in \{w_{q_1+q_2},w_n\}$.

We start with the assumption that $q_3\geq 4$.
By Remark~\ref{R:involution}, we know that 
$U_{I,J}^-$ is isomorphic to $U_{\theta(I),\theta(J)}^-$. 
Therefore, to prove that we can reduce to $q_3 \leq 3$,
we are going to work with the isomorphic poset 
$U_{\theta(I),\theta(J)}^-$, which is given by 
$Bl(P_{\theta(I)})=(p_2,3)$, $p_2\geq 3$ and 
$Bl(P_{\theta(J)})=(q_3,q_2,q_1)$, $q_1,q_2,q_3\geq 2$.
Note that $p_2=n-3$. 
Since $\theta(I)^c=\{s_{p_2}\}$, and 
$\theta(J)^c=\{s_{q_3},s_{q_3+q_2}\}$,
we see that if $w=w_1\dots w_n \in U_{\theta(I),\theta(J)}^-$, then 
\begin{enumerate}
\item for $i\in  [n-1]\setminus \{n-3\}$, $i$ comes before $i+1$ in $w$;
\item $w_1<\dots < w_{q_3}$, $w_{q_3+1}<\dots < w_{q_3+q_2}$,
$w_{q_3+q_2+1}<\dots < w_{n}$.
\end{enumerate}
This implies that $1\in \{w_1,w_{q_3+1},w_{q_3+q_2+1}\}$.
If $1$ appears as $w_{q_3+1}$ or 
$w_{q_3+q_2+1}$,
then we cannot fit $2,3,\dots, n-3$ in $w$ since they come
after 1 in $w$. Therefore, we have $w_1=1$. 
Then we remove 1 from all permutations in 
$U_{\theta(I),\theta(J)}^-$ and we reduce each remaining number 
by 1. This operation gives us a 
poset $U_{\theta(I)',\theta(J)'}^{'-}$, isomorphic to $U_{\theta(I),\theta(J)}^-$,
where $Bl(P_{\theta(I)'})=(p_2-1,3)$, $p_2-1\geq 3$ 
and $Bl(P_{\theta(J)'})=(q_3-1,q_2,q_1)$,
$q_1,q_2,q_3-1\geq 2$. 
Therefore, we can assume that $q_3\leq 3$.

Let us proceed with the assumption that $q_1\geq 4$,
and let $w=w_1\dots w_n$ be an element from $U_{I,J}^-$.
By condition (i), we know that $5$ appears either in 
the first segment $w_1\dots w_{q_1}$, or in the second
segment $w_{q_1+1}\dots w_{q_1+q_2}$. 
If it appears in the first segment, then 4 has to precede 
5 otherwise it creates a descent which gives a contradiction. 
If $5$ appears in the second segment $w_{q_1+1}\dots w_{q_1+q_2}$,
then we must have $w_5=5$ by conditions (i) and (ii), and by our assumption
that $q_1+1 \geq 5$.
In this case, condition (ii) shows that 4 has to be equal to $w_4$. 
These arguments show that if $q_1\geq 4$, 
then 4 precedes 5 in every element $w\in U_{I,J}^-$.
Therefore, removing 4 from $w$ and reducing every number
bigger than 4 by 1 
give us a new poset $U_{I', J'}^{-}$, 
isomorphic to $U_{I,J}^-$,
where $Bl(P_{I'})=(3,p_2-1)$, $p_2-1\geq 3$ 
and $Bl(P_{J'})=(q_1-1,q_2,q_3)$,
$q_1-1,q_2,q_3\geq 2$.

Now we assume that $q_2\geq 4$ along with $2\leq q_1,q_3\leq 3$. 
We will look for where in $w=w_1\dots w_n \in U_{I,J}^-$ the numbers 
$n-q_3$ and $n-q_3+1$ appear. Since $q_2\geq 4$, 
we see from conditions (i) and (ii) that $n-q_3$ appears in the segment 
$w_{q_1+1}<\dots < w_{q_1+q_2}$. 
We claim that if $w_k = w_{n-q_3}$ for some 
$k\in \{q_1+1,\dots, q_1+q_2\}$, then $w_{k+1} = w_{n-q_3+1}$.
This is clearly true if $n-q_3$ appears in the same segment 
$w_{q_1+1}\dots w_{q_1+q_2}$ since there is no descents 
within this segment. On the other hand, if $n-q_3+1$ appears 
in the segment $w_{q_1+q_2+1}<\dots < w_n$, then 
we must have $w_{q_1+q_2+1}=w_{n-q_3+1}=n-q_3+1$. 
But in this case, $w_{q_1+q_2+i}=n-q_3+i$, therefore, 
$w_{q_1+q_2} < w_{q_1+q_2+1}$. 
This implies that $n-q_3$ appears as the last 
entry $w_{q_1+q_2}$ of the segment $w_{q_1+1}\dots w_{q_1+q_2}$,
hence the proof of our claim follows. 
Now we know that $n-q_3$ and $n-q_3+1$ appear
in any $w\in U_{I,J}^-$ consecutively. Therefore, the removal of $n-q_3$ 
from $w$, and the reduction of all entries bigger than $n-q_3$ in $w$ by 1
gives a permutation in $S_{n-1}$.
Furthermore, this operation preserves  
the relative ordering (in Bruhat-Chevalley order) of the elements of $U_{I,J}^-$. 
In other words, we obtain a new poset $U_{I',J'}^{'-}$, 
isomorphic to $U_{I,J}^-$,
where $Bl(P_{I'})=(3,p_2-1)$, $p_2-1\geq 3$ 
and $Bl(P_{J'})=(q_1,q_2-1,q_3)$,
$q_1,q_2-1,q_3\geq 2$. 
These reduction arguments show that it suffices to 
consider the following eight cases only:

\begin{enumerate}
\item $Bl(P_I)=(3,3)$, $Bl(P_J)=(2,2,2)$;
\item $Bl(P_I)=(3,4)$, $Bl(P_J)=(2,2,3)$;
\item $Bl(P_I)=(3,4)$, $Bl(P_J)=(3,2,2)$; 
\item $Bl(P_I)=(3,4)$, $Bl(P_J)=(2,3,2)$;
\item $Bl(P_I)=(3,5)$, $Bl(P_J)=(2,3,3)$;
\item $Bl(P_I)=(3,5)$, $Bl(P_J)=(3,2,3)$;
\item $Bl(P_I)=(3,5)$, $Bl(P_J)=(3,3,2)$;
\item $Bl(P_I)=(3,6)$, $Bl(P_J)=(3,3,3)$.
\end{enumerate}
The Hasse diagrams of the posets corresponding to these eight cases 
are given by the diagrams $P.1 - P.8$ in Figure~\ref{F:28}.

\subsection{$Bl(P_I)=(p_1,p_2)$, $p_1,p_2\geq 3$ and $Bl(P_J)=(2,2,q_3)$, $q_3\geq 2$.}

First, we assume that $p_2\geq 5$,
and we apply $\theta$ to $I$ and $J$. 
Then $\theta(I)^{c}=\{s_{p_2}\}$, 
and $\theta(J)^{c}=\{s_{n-4},s_{n-2}\}$,
we see that if $w=w_1\dots w_n \in U_{\theta(I),\theta(J)}^-$, then 
\begin{enumerate}
\item[(i)] for $i\in  [n-1]\setminus \{p_2\}$, 
$i$ comes before $i+1$ in $w$;
\item[(ii)] $w_1<\dots < w_{n-4}$, $w_{n-3}< w_{n-2}$,
$w_{n-1}<w_{n}$.
\end{enumerate}
This means that $1$ is contained in $\{w_1,w_{n-3},w_{n-1}\}$.
Recall that $p_2\geq 5$. Thus, we cannot place 
the sequence $1,2,\dots, p_2$ in $w$ as an increasing substring 
unless $w_1=1$. So, $w$ starts with 1. Since this is true for all 
elements of $U_{\theta(I),\theta(J)}^-$, by first removing $w_1= 1$ from 
all $w\in U_{\theta(I),\theta(J)}^-$, and then reducing the remaining 
entries by 1, we obtain an isomorphic 
poset $U_{\theta(I)',\theta(J)'}^-$ in $S_{n-1}$,
where $\theta(I)^{'c}=\{ s_{p_2-1} \}$ and 
$\theta(J)^{'c}=\{s_{n-4},s_{n-2}\}$.
Therefore, we see that we can assume $p_2\leq 4$.

We now proceed with the assumption that $p_1 \geq 5$ and that $p_2\leq 4$.
If $w=w_1\dots w_n \in U_{I,J}^-$, then 
\begin{enumerate}
\item for $i\in  [n-1]\setminus \{p_1\}$, 
$i$ comes before $i+1$ in $w$;
\item $w_1<w_2,\ w_3<w_4,\ w_5<\dots < w_n$.
\end{enumerate}
We will look for where in $w=w_1\dots w_n$ the numbers 
$p_1-1$ and $p_1$ appear. Since $p_1\geq 5$, 
we see from conditions 1 and 2 that $p_1$ appears in the segment 
$w_5<w_6<\dots < w_n$. If $w_k = p_1$ and $k>5$,
then clearly $w_{k-1}=p_1-1$ otherwise we must have a descent in
the segment $w_5w_6\dots w_n$, which would contradict with Condition 2. 
On the other hand, if $w_5=p_1$, then we see that $5=p_1$,
hence $w_4=p_1-1$. In both of these cases, we see that 
if $w_k=p_1$, then $w_{k-1}=p_1-1$. 
Now, by removing $p_1$ from $w\in U_{I,J}^-$ 
and reducing by 1 all entries $w_j$ with $w_j > p_1$,
we obtain a poset $U_{I',J'}^-$,
isomorphic to $U_{I,J}^-$, in $S_{n-1}$. 
Furthermore, $Bl(P_{I'})=(p_1-1,p_2)$, $p_1-1,p_2\geq 3$ 
and $Bl(P_{J'})=(2,2,q_3-1)$, $q_3-1\geq 2$.
In other words, we can assume that $p_1\leq 4$.

These two reduction arguments show that it suffices to 
consider the following four cases only: 
\begin{enumerate}
\item $Bl(P_I)=(3,3)$, $Bl(P_J)=(2,2,2)$;
\item $Bl(P_I)=(3,4)$, $Bl(P_J)=(2,2,3)$;
\item $Bl(P_I)=(4,3)$, $Bl(P_J)=(2,2,3)$;
\item $Bl(P_I)=(4,4)$, $Bl(P_J)=(2,2,4)$.
\end{enumerate}
The Hasse diagrams of the posets corresponding to these four 
cases are given by the diagrams $P.1, P.2, P.3$, and $P.6$ of Figure~\ref{F:28}.

\subsection{$Bl(P_I)=(p_1,p_2)$, $p_1,p_2\geq 3$ and $Bl(P_J)=(2,q_2,2)$, $q_2\geq 2$.}

First, we assume that $p_1\geq 5$.
Since $I^c=\{s_{p_1}\}$, $J^c=\{s_{2},s_{n-2}\}$ in $R_{n-1}$, 
we see that if $w=w_1\dots w_n \in U_{I,J}^-$, then 
\begin{enumerate}
\item[(i)] for $i\in [n-1]\setminus \{p_1\}$, 
$i$ comes before $i+1$ in $w$;
\item[(ii)] $w_1<w_2,\ w_3< \dots < w_{n-2}$,\ $w_{n-1}< w_{n}$.
\end{enumerate}

We look for the positions of $p_1-3$ and $p_1-2$. 
Since $p_1\geq 5$, we see from condition (i) that $p_1-2$ 
appears in the segment $w_3w_4\dots w_{n-2}$.
If $w_k = p_1-2$ for some $k>3$, then we see 
that $p_1-3$ must also be in the same segment, 
hence, we must have that $w_{k-1}=p_1-3$. 
If $w_3=p_1-2$, then, by conditions (i) and (ii), 
we have only one choice that $p_1=5$, and 
$p_1-3=2=w_2$. 
In both of these two cases we see that $p_1-3$ must 
come immediately before $p_1-2$ in every $w\in U_{I,J}^-$.
Therefore, by removing $p_1-2$ from $w$ and 
reducing every entry which is greater than $p_1-2$ 
by 1, we do not change the structure of the underlying poset; 
we obtain a poset $U_{I',J'}^-$ in $S_{n-1}$ 
such that $Bl(P_{I'})=(p_1-1,p_2)$, $p_1-1,p_2\geq 3$ 
and $Bl(P_{J'})=(2,q_2-1,2)$, $q_2-1\geq 2$.
In other words, we can assume that $p_1 \leq 4$.

For $p_2\geq 5$, we repeat the same arguments after 
applying $\theta$ to $I$ and $J$. 
Therefore, without loss of generality we can assume 
that $3\leq p_1,p_2\leq 4$. 
This reduction argument shows that our poset is isomorphic 
to one of the following three cases:
\begin{enumerate}
\item $Bl(P_I)=(3,3)$, $Bl(P_J)=(2,2,2)$;
\item $Bl(P_I)=(3,4)$, $Bl(P_J)=(2,3,2)$;
\item $Bl(P_I)=(4,4)$, $Bl(P_J)=(2,4,2)$.
\end{enumerate}
The Hasse diagrams of the posets corresponding to these three 
cases are given by the diagrams $P.1, P.4$ and $P.9$ in Figure~\ref{F:28}.

\subsection{$Bl(P_I)=(2,p_2)$, $p_2\geq 3$ and $Bl(P_J)=(q_1,q_2,q_3,q_4)$.}

Let us first assume that $q_4 \geq 3$. 
Since $I^c=\{s_2\}$, $J^c=\{s_{q_1},s_{q_1+q_2},s_{q_1+q_2+q_3}\}$ in $R_{n-1}$, 
we see that if $w=w_1\dots w_n \in U_{I,J}^-$, then 
\begin{enumerate}
\item[(i)] for $i\in  [n-1]\setminus \{2\}$, $i$ comes before $i+1$ in $w$;
\item[(ii)] $w_1< \dots < w_{q_1}$, $w_{q_1+1}< \dots < w_{q_1+q_2}$, 
$w_{q_1+q_2+1}< \dots < w_{q_1+q_2+q_3}$, and 
$w_{q_1+q_2+q_3+1}< \dots < w_{n}$.
\end{enumerate}
This implies that $n\in \{w_{q_1},w_{q_1+q_2},w_{q_1+q_2+q_3},w_n\}$.
By (i) we know that $n$ is preceded by $3,\dots, n-1$, which prevents the 
possibilities $n\in \{w_{q_1},w_{q_1+q_2},w_{q_1+q_2+q_3}\}$. 
Therefore, $w_n=n$. 
Thus, by removing $n$ from $w\in U_{I,J}^-$, 
we do not change the structure of the underlying poset; 
we obtain a poset $U_{I',J'}^-$ in $S_{n-1}$, 
which is isomorphic to $U_{I,J}^-$, 
such that $Bl(P_{I'})=(2,p_2-1)$, $p_2-1\geq 3$ 
and $Bl(P_{J'})=(q_1,q_2,q_3,q_4-1)$. In other words, 
we can assume without loss of generality that $1\leq q_4\leq 2$. 

We proceed with the assumption that $q_3\geq 3$. 
Then we look at the relative positions of the numbers $m:=q_1+q_2+q_3$ 
and $m+1$ in $w$. 
Since we assumed that $1\leq q_4\leq 2$, we have 
$n\in \{ w_{m+1}, w_n\}$. 
If $n=w_{m+1}$, then the following implication is obvious:
$$
w_k=m \implies w_{k+1}=m+1.
$$
On the other hand, if $n=w_n$, then since $q_3\geq 3$,
we know that $m+1$ has to appear in the 
following segment of $w$: 
$w_{q_1+q_2+1} \dots  w_{q_1+q_2+q_3}$.
In particular, we have one of the following cases: 
$$
w_{q_1+q_2+q_3-i}=m\ \text{ and }\ w_{q_1+q_2+q_3-i+1}=m+1
$$
for $i=0,1$. 
Therefore, $m$ and $m+1$ appear as consecutive terms in $w$,
furthermore, $m$ appears in $w_{q_1+q_2+1} \dots  w_{q_1+q_2+q_3}$.
In this case, by removing $m$ from $w$ and reducing 
every number greater than $m$ by 1, we obtain a poset 
$U_{I',J'}^-$ in $S_{n-1}$, which is isomorphic to $U_{I,J}^-$, 
such that $Bl(P_{I'})=(2,p_2-1)$, $p_2-1\geq 3$ 
and $Bl(P_{J'})=(q_1,q_2,q_3-1,q_4)$. In other words, 
we can assume without loss of generality that $1\leq q_3\leq 2$ as well.

Next, we proceed with the assumptions that $q_2\geq 3$ and $1\leq q_3,q_4 \leq 2$. 
In this case, after applying the involution $\theta$ to $I$ and $J$, we assume 
that $Bl(P_I)=(p_2,2)$, $p_2\geq 3$ and $Bl(P_J)=(q_4,q_3,q_2,q_1)$, where 
$q_2\geq 3$ and $1\leq q_3,q_4 \leq 2$. 
In other words, we have one of the following four possibilities for the first few terms
of $J$:
\begin{enumerate}
\item $s_1,s_3,s_5,s_6 \in J$ and $s_2,s_4\notin J$, or
\item $s_1,s_4,s_5 \in J$ and $s_2,s_3\notin J$, or 
\item $s_2,s_4,s_5 \in J$ and $s_1,s_3\notin J$, or
\item $s_3,s_4 \in J$ and $s_1,s_2\notin J$.
\end{enumerate}
In the first case, we have that 
$$
w_k=4 \implies w_{k+1}=5
$$
for some $k\geq 1$. 
In the second case, we have 
$$
w_k=3 \implies w_{k+1}=4
$$
for some $k\geq 1$. 
In the third case, we have 
$$
w_k=3 \implies w_{k+1}=4
$$
for some $k\geq 1$. 
Finally, in the fourth case, we have
$$
w_k=2 \implies w_{k+1}=3
$$
for some $k\geq 1$. 
In all of these cases, removing $w_{k+1}$ from $w$ and 
reducing every number that is greater than $w_{k+1}$ by 1
give a poset 
$U_{I',J'}^-$ in $S_{n-1}$, which is isomorphic to $U_{I,J}^-$, 
such that $Bl(P_{I'})=(p_2-1,2)$, $p_2-1\geq 3$ and $Bl(P_{J'})=(q_4,q_3,q_2-1,q_1)$. 
In other words, we can assume without loss of generality that $1\leq q_2 \leq 2$.

Finally, we assume that $q_1\geq 3$ and $1\leq q_2,q_3,q_4\leq 2$. 
The proof of this case develops similar to the previous case;
we apply $\theta$ to $I$ and $J$; we assume 
that $Bl(P_I)=(p_2,2)$, $p_2\geq 3$ and $Bl(P_J)=(q_4,q_3,q_2,q_1)$, 
where $q_1\geq 3$ and $1\leq q_2,q_3,q_4 \leq 2$. 
This time we have 8 possibilities, instead of 4 as in the previous case.
In each of these eight cases, we consider the simple reflection $s_j$ with 
smallest index $j$ among the elements of $J$ associated to its block of size $q_1$. 
Then, as in the previous case, 
$$
w_k=j-1 \implies w_{k+1}=j
$$
for some $k\geq 1$. 
Therefore, removing $j$ from $w$ and 
reducing every number that is greater than $j$ by 1
give a poset 
$U_{I',J'}^-$ in $S_{n-1}$, isomorphic to $U_{I,J}^-$, 
such that $Bl(P_{I'})=(p_2-1,2)$, $p_2-1\geq 3$ 
and $Bl(P_{J'})=(q_4,q_3,q_2,q_1-1)$. In other words, 
we can assume without loss of generality that $1\leq q_1 \leq 2$.
Now we know that if 
$Bl(P_I)=(2,p_2)$, $p_2\geq 3$ and $Bl(P_J)=(q_1,q_2,q_3,q_4)$, then
$U_{I,J}^-$ is isomorphic to one of the following 15 cases:
\begin{enumerate}
\item $Bl(P_I)=(2,3)$ and $Bl(P_J)=(1,1,1,2)$;
\item $Bl(P_I)=(2,3)$ and $Bl(P_J)=(1,1,2,1)$;
\item $Bl(P_I)=(2,3)$ and $Bl(P_J)=(1,2,1,1)$;
\item $Bl(P_I)=(2,3)$ and $Bl(P_J)=(2,1,1,1)$;
\item $Bl(P_I)=(2,4)$ and $Bl(P_J)=(1,1,2,2)$;
\item $Bl(P_I)=(2,4)$ and $Bl(P_J)=(1,2,1,2)$;
\item $Bl(P_I)=(2,4)$ and $Bl(P_J)=(2,1,1,2)$;
\item $Bl(P_I)=(2,4)$ and $Bl(P_J)=(1,2,2,1)$;
\item $Bl(P_I)=(2,4)$ and $Bl(P_J)=(2,1,2,1)$;
\item $Bl(P_I)=(2,4)$ and $Bl(P_J)=(2,2,1,1)$;
\item $Bl(P_I)=(2,5)$ and $Bl(P_J)=(1,2,2,2)$;
\item $Bl(P_I)=(2,5)$ and $Bl(P_J)=(2,1,2,2)$;
\item $Bl(P_I)=(2,5)$ and $Bl(P_J)=(2,2,1,2)$;
\item $Bl(P_I)=(2,5)$ and $Bl(P_J)=(2,2,2,1)$;
\item $Bl(P_I)=(2,6)$ and $Bl(P_J)=(2,2,2,2)$.
\end{enumerate}
The Hasse diagrams of the posets corresponding to these 15 cases 
are given by the diagrams 
$P.10, P.11, P.12, P.13, P.14, P.4, P.15, P.4, P.15, P.14,
P.5, P.16, P.17, P.7, P.8$ in Figure~\ref{F:28}.
Note that several Hasse diagrams appear multiple times in this list.

\subsection{$Bl(P_I)=(p_1,p_2)$, $p_1,p_2\geq 2$ and $Bl(P_J)=(1,1,1,q_4)$.} 

We consider this situation in two different cases: 
\begin{enumerate}
\item[(a)] $Bl(P_I)=(2,2)$ and $Bl(P_J)=(1,1,1,1)$;
\item[(b)] $Bl(P_I)=(p_1,p_2)$, $p_1,p_2\geq 2$ 
and $Bl(P_J)=(1,1,1,q_4)$, $q_4 \geq 2$.
\end{enumerate}

We explain the reduction argument for (b);
we claim that we can assume $2\leq p_1,p_2 \leq 3$.

First, we assume that $p_2\geq 4$.
Since $I^c=\{s_{p_1}\}$, $J^c=\{s_1,s_2,s_3\}$ in $R_{n-1}$, 
we see that if $w=w_1\dots w_n \in U_{I,J}^-$, then 
\begin{enumerate}
\item[(i)] for $i\in  [n-1]\setminus \{p_1\}$, 
$i$ comes before $i+1$ in $w$;
\item[(ii)] $w_4< \dots <  w_{n}$.
\end{enumerate}
Therefore, $n\in \{w_1,w_2,w_3,w_n\}$. But there are at least $p_2-1 \geq 3$
numbers before $n$ in $w$, therefore, $n$ cannot appear 
in $\{w_1,w_2,w_3\}$. This means that $n$ is equal to $w_n$. 
Now we see that removing $n$ from $w$, for all $w\in U_{I,J}^-$
gives us an isomorphic poset $U_{I',J'}^-$, 
where $Bl(P_{I'})=(p_1,p_2-1)$, $p_2-1,p_1\geq 2$ 
and $Bl(P_{J'})=(1,1,1,q_4-1)$, $q_4-1\geq 2$.

We now proceed with the assumption that $p_1\geq 4$. 
In this case, we look at the relative positions of numbers 3 and 4.
If $3$ appears in the segment $w_4 w_5 \dots w_n$, 
then 3 is immediately followed by 4 since there are no descents in this portion of $w$.
On the other hand, if 3 does not appear in the segment $w_4 w_5 \dots w_n$,
then it can only appear at $w_3$ since in this case it has to be preceded by 1 and 2
by condition (i). But then, 4 has to appear as $w_4$, otherwise, there would 
be a descent in $w_4 w_5 \dots w_n$. This argument shows that the 
numbers 3 and 4 appear in $w$ consecutively. Hence, if we remove 4 from $w$,
and reduce every number greater than 4 by 1, 
then we do not change the Bruhat-Chevalley order. In other words, 
we obtain a poset $U_{I',J'}^-$, isomorphic to $U_{I,J}^-$, 
where $Bl(P_{I'})=(p_1-1,p_2)$, $p_1-1,p_2\geq 2$ 
and $Bl(P_{J'})=(1,1,1,q_4-1)$, $q_4-1\geq 2$. Therefore, we can 
assume that $p_1\leq 3$. 
Consequently, we see in this case that there are only the following 
four possibilities:

\begin{enumerate}
\item $Bl(P_I)=(2,2)$ and $Bl(P_J)=(1,1,1,1)$;
\item $Bl(P_I)=(2,3)$ and $Bl(P_J)=(1,1,1,2)$; 
\item $Bl(P_I)=(3,2)$ and $Bl(P_J)=(1,1,1,2)$; 
\item $Bl(P_I)=(3,3)$ and $Bl(P_J)=(1,1,1,3)$.
\end{enumerate}
The Hasse diagrams of the corresponding posets of these four cases 
are given by the diagrams $P.18, P.10, P.13$, and $P.19$ in Figure~\ref{F:28}.

\subsection{$Bl(P_I)=(p_1,p_2)$, $p_1,p_2\geq 2$ and $Bl(P_J)=(1,1,q_3,1)$, $q_3\geq 2$.} 

By arguing as in the previous cases,
we see that each subcase reduces to the one of the following three subcases: 

\begin{enumerate}
\item $Bl(P_I)=(2,3)$ and $Bl(P_J)=(1,1,2,1)$;
\item $Bl(P_I)=(3,2)$ and $Bl(P_J)=(1,1,2,1)$;
\item $Bl(P_I)=(3,3)$ and $Bl(P_J)=(1,1,3,1)$.
\end{enumerate}
The Hasse diagrams of the corresponding posets of these three cases 
are given by $P.11, P.12$ and $P.20$ in Figure~\ref{F:28}.

\subsection{$Bl(P_I)=(1,p_2,1)$, $p_2\geq 2$ and $Bl(P_J)=(q_1,q_2,q_3)$.}

In this case, by the appropriate reduction arguments as in the previous cases
we may assume that $q_1,q_2,q_3\leq 2$. 
Then we obtain the following five cases: 

\begin{enumerate}
\item $Bl(P_I)=(1,2,1)$ and $Bl(P_J)=(1,1,2)$;
\item $Bl(P_I)=(1,2,1)$ and $Bl(P_J)=(1,2,1)$;
\item $Bl(P_I)=(1,3,1)$ and $Bl(P_J)=(1,2,2)$;
\item $Bl(P_I)=(1,3,1)$ and $Bl(P_J)=(2,1,2)$;
\item $Bl(P_I)=(1,4,1)$ and $Bl(P_J)=(2,2,2)$.
\end{enumerate}
The Hasse diagrams of the corresponding posets of these five cases 
are given by $P.21$, $P.1$, $P.4$, $P.22$, and $P.9$ in Figure~\ref{F:28}.

\subsection{$Bl(P_I)=(1,1,p_3)$, $p_3\geq 2$ and $Bl(P_J)=(q_1,q_2,q_3)$.}

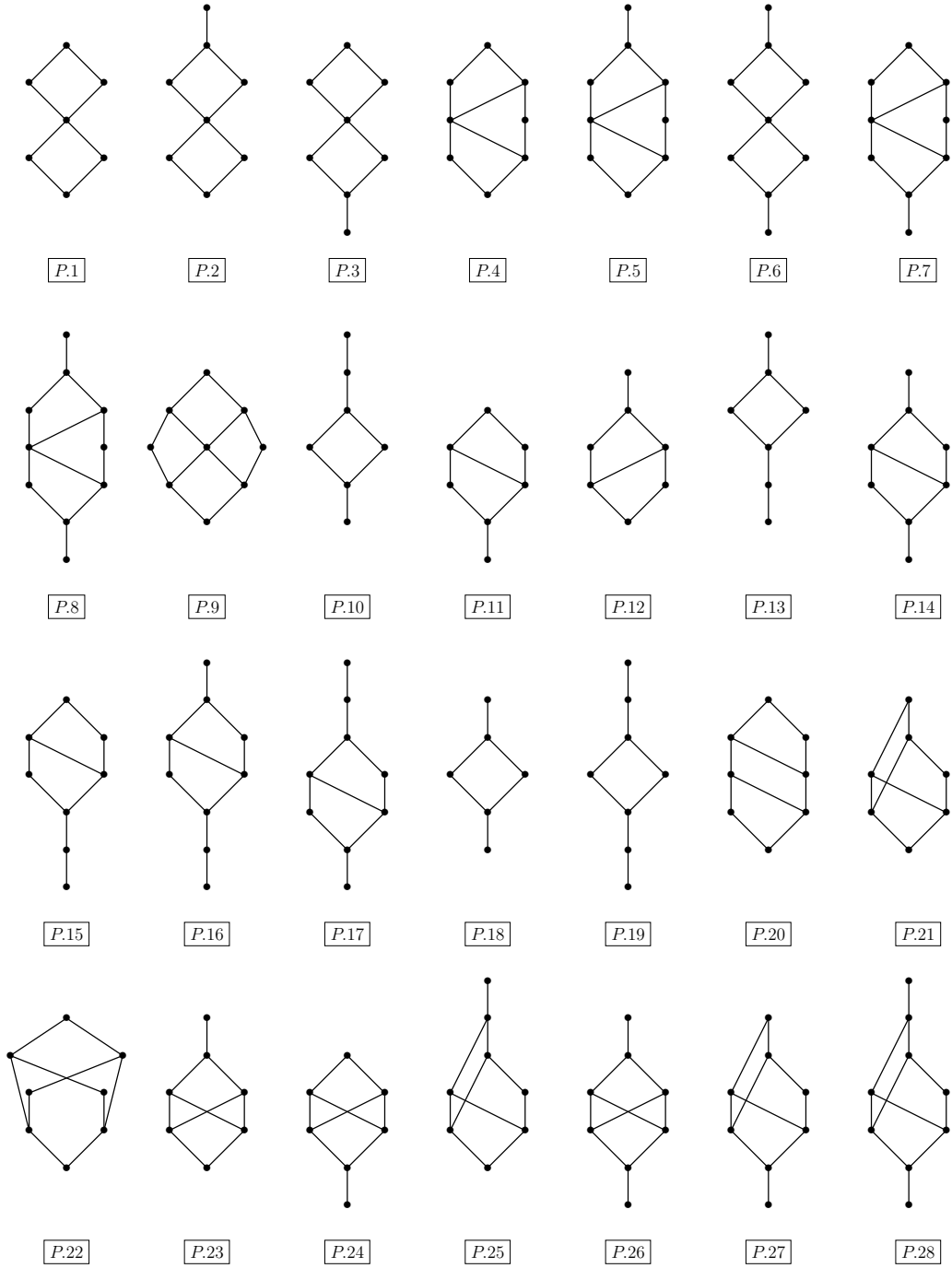
\begin{figure}[htp]

\centering
\scalebox{.6}{
\begin{tikzpicture}[scale=.45]

\draw (-22.5,46) node [draw] {$P.1$};
\draw (-15,46) node [draw] {$P.2$};
\draw (-7.5,46) node [draw] {$P.3$};
\draw (0,46) node [draw] {$P.4$};
\draw (7.5,46) node [draw] {$P.5$};
\draw (15,46) node [draw] {$P.6$};
\draw (23,46) node [draw] {$P.7$};

\draw (-22.5,28) node [draw] {$P.8$};
\draw (-15,28) node [draw] {$P.9$};
\draw (-7.5,28) node [draw] {$P.10$};
\draw (0,28) node [draw] {$P.11$};
\draw (7.5,28) node [draw] {$P.12$};
\draw (15,28) node [draw] {$P.13$};
\draw (23,28) node [draw] {$P.14$};

\draw (-22.5,10.5) node [draw] {$P.15$};
\draw (-15,10.5) node [draw] {$P.16$};
\draw (-7.5,10.5) node [draw] {$P.17$};
\draw (0,10.5) node [draw] {$P.18$};
\draw (7.5,10.5) node [draw] {$P.19$};
\draw (15,10.5) node [draw] {$P.20$};
\draw (23,10.5) node [draw] {$P.21$};

\draw (-22.5,-6.5) node [draw] {$P.22$};
\draw (-15,-6.5) node [draw] {$P.23$};
\draw (-7.5,-6.5) node [draw] {$P.24$};
\draw (0,-6.5) node [draw] {$P.25$};
\draw (7.5,-6.5) node [draw] {$P.26$};
\draw (15,-6.5) node [draw] {$P.27$};
\draw (23,-6.5) node [draw] {$P.28$};

\begin{scope}[xshift=-22.5cm, yshift=50cm]
\draw[-,thick] (0,0) to (-2,2) to (0,4) to (-2,6) to (0,8) to (2,6) to (0,4) to (2,2) to (0,0);
\node at (0,0) {$\bullet$};
\node at (-2,2) {$\bullet$};
\node at (0,4) {$\bullet$};
\node at (2,2) {$\bullet$};
\node at (2,6) {$\bullet$};
\node at (-2,6) {$\bullet$};
\node at (0,8) {$\bullet$};
\end{scope}

\begin{scope}[xshift=-15cm, yshift=50cm]
\draw[-,thick] (0,0) to (-2,2) to (0,4) to (-2,6) to (0,8) to (2,6) to (0,4) to (2,2) to (0,0);
\draw[-,thick] (0,8) to (0,10) ;
\node at (0,0) {$\bullet$};
\node at (-2,2) {$\bullet$};
\node at (0,4) {$\bullet$};
\node at (2,2) {$\bullet$};
\node at (2,6) {$\bullet$};
\node at (-2,6) {$\bullet$};
\node at (0,8) {$\bullet$};
\node at (0,10) {$\bullet$};
\end{scope}

\begin{scope}[xshift=-7.5cm, yshift=50cm]
\draw[-,thick] (0,0) to (-2,2) to (0,4) to (-2,6) to (0,8) to (2,6) to (0,4) to (2,2) to (0,0);
\draw[-,thick] (0,0) to (0,-2);
\node at (0,0) {$\bullet$};
\node at (-2,2) {$\bullet$};
\node at (0,4) {$\bullet$};
\node at (2,2) {$\bullet$};
\node at (2,6) {$\bullet$};
\node at (-2,6) {$\bullet$};
\node at (0,8) {$\bullet$};
\node at (0,-2) {$\bullet$};
\end{scope}

\begin{scope}[xshift=0cm, yshift=50cm]

\draw[-,thick] (0,0) to (-2,2) to (-2,4) to (-2,6) to (0,8) to (2,6) to (2,4) to (2,2) to (0,0);
\draw[-,thick] (2,2) to (-2,4) to (2,6) ;
\node at (0,0) {$\bullet$};
\node at (-2,2) {$\bullet$};

\node at (2,2) {$\bullet$};
\node at (2,4) {$\bullet$};
\node at (-2,4) {$\bullet$};
\node at (2,6) {$\bullet$};
\node at (-2,6) {$\bullet$};
\node at (0,8) {$\bullet$};

\end{scope}

\begin{scope}[xshift=7.5cm, yshift=50cm]

\draw[-,thick] (0,0) to (-2,2) to (-2,4) to (-2,6) to (0,8) to (2,6) to (2,4) to (2,2) to (0,0);
\draw[-,thick] (2,2) to (-2,4) to (2,6) ;
\draw[-,thick] (0,8) to (0,10) ;
\node at (0,0) {$\bullet$};
\node at (-2,2) {$\bullet$};

\node at (2,2) {$\bullet$};
\node at (2,4) {$\bullet$};
\node at (-2,4) {$\bullet$};
\node at (2,6) {$\bullet$};
\node at (-2,6) {$\bullet$};
\node at (0,8) {$\bullet$};
\node at (0,10) {$\bullet$};
\end{scope}

\begin{scope}[xshift=15cm, yshift=50cm]
\draw[-,thick] (0,0) to (-2,2) to (0,4) to (-2,6) to (0,8) to (2,6) to (0,4) to (2,2) to (0,0);
\draw[-,thick] (0,0) to (0,-2);
\draw[-,thick] (0,8) to (0,10);
\node at (0,0) {$\bullet$};
\node at (-2,2) {$\bullet$};
\node at (0,4) {$\bullet$};
\node at (2,2) {$\bullet$};
\node at (2,6) {$\bullet$};
\node at (-2,6) {$\bullet$};
\node at (0,8) {$\bullet$};
\node at (0,-2) {$\bullet$};
\node at (0,10) {$\bullet$};
\end{scope}

\begin{scope}[xshift=22.5cm, yshift=50cm]

\draw[-,thick] (0,0) to (-2,2) to (-2,4) to (-2,6) to (0,8) to (2,6) to (2,4) to (2,2) to (0,0);
\draw[-,thick] (2,2) to (-2,4) to (2,6) ;
\draw[-,thick] (0,0) to (0,-2) ;
\node at (0,0) {$\bullet$};
\node at (-2,2) {$\bullet$};

\node at (2,2) {$\bullet$};
\node at (2,4) {$\bullet$};
\node at (-2,4) {$\bullet$};
\node at (2,6) {$\bullet$};
\node at (-2,6) {$\bullet$};
\node at (0,8) {$\bullet$};
\node at (0,-2) {$\bullet$};
\end{scope}


\begin{scope}[xshift=-22.5cm, yshift=32.5cm]
\draw[-,thick] (0,0) to (-2,2) to (-2,4) to (-2,6) to (0,8) to (2,6) to (2,4) to (2,2) to (0,0);
\draw[-,thick] (2,2) to (-2,4) to (2,6) ;
\draw[-,thick] (0,8) to (0,10) ;
\draw[-,thick] (0,0) to (0,-2) ;
\node at (0,0) {$\bullet$};
\node at (-2,2) {$\bullet$};
\node at (2,2) {$\bullet$};
\node at (2,4) {$\bullet$};
\node at (-2,4) {$\bullet$};
\node at (2,6) {$\bullet$};
\node at (-2,6) {$\bullet$};
\node at (0,8) {$\bullet$};
\node at (0,10) {$\bullet$};
\node at (0,-2) {$\bullet$};
\end{scope}

\begin{scope}[xshift=-15cm, yshift=32.5cm]
\draw[-,thick] (0,0) to (-2,2) to (-3,4) to (-2,6) to (0,8) to (2,6) to (3,4) to (2,2) to (0,0);
\draw[-,thick] (-2,2) to (0,4) to (2,6) ;
\draw[-,thick] (2,2) to (0,4) to (-2,6) ;
\node at (0,0) {$\bullet$};
\node at (-2,2) {$\bullet$};
\node at (2,2) {$\bullet$};
\node at (-3,4) {$\bullet$};
\node at (0,4) {$\bullet$};
\node at (3,4) {$\bullet$};
\node at (2,6) {$\bullet$};
\node at (-2,6) {$\bullet$};
\node at (0,8) {$\bullet$};
\end{scope}

\begin{scope}[xshift=-7.5cm, yshift=32.5cm]
\draw[-,thick] (0,0) to (0,2) to (2,4) to (0,6) to (0,8) to (0,10);
\draw[-,thick] (0,2) to (-2,4) to (0,6);
\node at (0,0) {$\bullet$};
\node at (0,2) {$\bullet$};
\node at (2,4) {$\bullet$};
\node at (0,6) {$\bullet$};
\node at (-2,4) {$\bullet$};
\node at (0,8) {$\bullet$};
\node at (0,10) {$\bullet$};
\end{scope}

\begin{scope}[xshift=0cm, yshift=32.5cm]
\draw[-,thick] (0,0) to (-2,2) to (-2,4) to (0,6) to (2,4) to (2,2) to (0,0);
\draw[-,thick] (2,2) to (-2,4);
\draw[-,thick] (0,0) to (0,-2);
\node at (0,0) {$\bullet$};
\node at (2,2) {$\bullet$};
\node at (-2,2) {$\bullet$};
\node at (2,4) {$\bullet$};
\node at (-2,4) {$\bullet$};
\node at (0,6) {$\bullet$};
\node at (0,-2) {$\bullet$};
\end{scope}

\begin{scope}[xshift=7.5cm, yshift=32.5cm]
\draw[-,thick] (0,0) to (-2,2) to (-2,4) to (0,6) to (2,4) to (2,2) to (0,0);
\draw[-,thick] (-2,2) to (2,4);
\draw[-,thick] (0,6) to (0,8);
\node at (0,0) {$\bullet$};
\node at (2,2) {$\bullet$};
\node at (-2,2) {$\bullet$};
\node at (2,4) {$\bullet$};
\node at (-2,4) {$\bullet$};
\node at (0,6) {$\bullet$};
\node at (0,8) {$\bullet$};
\end{scope}

\begin{scope}[xshift=15cm, yshift=32.5cm]
\draw[-,thick] (0,0) to (0,2) to (0,4) to (2,6) to (0,8) to (0,10);
\draw[-,thick] (0,4) to (-2,6) to (0,8);
\node at (0,0) {$\bullet$};
\node at (0,2) {$\bullet$};
\node at (0,4) {$\bullet$};
\node at (2,6) {$\bullet$};
\node at (-2,6) {$\bullet$};
\node at (0,8) {$\bullet$};
\node at (0,10) {$\bullet$};
\end{scope}

\begin{scope}[xshift=22.5cm, yshift=32.5cm]
\draw[-,thick] (0,0) to (-2,2) to (-2,4) to (0,6) to (2,4) to (2,2) to (0,0);
\draw[-,thick] (2,2) to (-2,4);
\draw[-,thick] (0,0) to (0,-2);
\draw[-,thick] (0,6) to (0,8);
\node at (0,0) {$\bullet$};
\node at (2,2) {$\bullet$};
\node at (-2,2) {$\bullet$};
\node at (2,4) {$\bullet$};
\node at (-2,4) {$\bullet$};
\node at (0,6) {$\bullet$};
\node at (0,8) {$\bullet$};
\node at (0,-2) {$\bullet$};
\end{scope}


\begin{scope}[xshift=-22.5cm, yshift=15cm]
\draw[-,thick] (0,-2) to (0,0) to (0,2) to (-2,4) to (-2,6) to (0,8) to (2,6) to (2,4) to (0,2);
\draw[-,thick] (2,4) to (-2,6);
\node at (0,8) {$\bullet$};
\node at (2,6) {$\bullet$};
\node at (-2,6) {$\bullet$};
\node at (2,4) {$\bullet$};
\node at (-2,4) {$\bullet$};
\node at (0,2) {$\bullet$};
\node at (0,0) {$\bullet$};
\node at (0,-2) {$\bullet$};
\end{scope}

\begin{scope}[xshift=-15cm, yshift=15cm]
\draw[-,thick] (0,-2) to (0,0) to (0,2) to (-2,4) to (-2,6) to (0,8) to (2,6) to (2,4) to (0,2);
\draw[-,thick] (0,8) to (0,10);
\draw[-,thick] (2,4) to (-2,6);
\node at (0,10) {$\bullet$};
\node at (0,8) {$\bullet$};
\node at (2,6) {$\bullet$};
\node at (-2,6) {$\bullet$};
\node at (2,4) {$\bullet$};
\node at (-2,4) {$\bullet$};
\node at (0,2) {$\bullet$};
\node at (0,0) {$\bullet$};
\node at (0,-2) {$\bullet$};
\end{scope}

\begin{scope}[xshift=-7.5cm, yshift=15cm]
\draw[-,thick] (0,0) to (-2,2) to (-2,4) to (0,6) to (2,4) to (2,2) to (0,0);
\draw[-,thick] (2,2) to (-2,4);
\draw[-,thick] (0,0) to (0,-2);
\draw[-,thick] (0,6) to (0,8);
\draw[-,thick] (0,8) to (0,10);
\node at (0,0) {$\bullet$};
\node at (2,2) {$\bullet$};
\node at (-2,2) {$\bullet$};
\node at (2,4) {$\bullet$};
\node at (-2,4) {$\bullet$};
\node at (0,6) {$\bullet$};
\node at (0,8) {$\bullet$};
\node at (0,10) {$\bullet$};
\node at (0,-2) {$\bullet$};
\end{scope}

\begin{scope}[xshift=0cm, yshift=15cm]
\draw[-,thick] (0,0) to (0,2) to (2,4) to (0,6) to (0,8);
\draw[-,thick] (0,2) to (-2,4) to (0,6);
\node at (0,0) {$\bullet$};
\node at (0,2) {$\bullet$};
\node at (2,4) {$\bullet$};
\node at (0,6) {$\bullet$};
\node at (-2,4) {$\bullet$};
\node at (0,8) {$\bullet$};
\end{scope}

\begin{scope}[xshift=7.5cm, yshift=15cm]
\draw[-,thick] (0,0) to (0,2) to (2,4) to (0,6) to (0,8) to (0,10);
\draw[-,thick] (0,2) to (-2,4) to (0,6);
\draw[-,thick] (0,0) to (0,-2);
\node at (0,0) {$\bullet$};
\node at (0,2) {$\bullet$};
\node at (2,4) {$\bullet$};
\node at (-2,4) {$\bullet$};
\node at (0,6) {$\bullet$};
\node at (0,8) {$\bullet$};
\node at (0,10) {$\bullet$};
\node at (0,-2) {$\bullet$};
\end{scope}

\begin{scope}[xshift=15cm, yshift=15cm]
\draw[-,thick] (0,0) to (-2,2) to (-2,4) to (-2,6) to (0,8) to (2,6) to (2,4) to (2,2) to (0,0);
\draw[-,thick] (2,2) to (-2,4);
\draw[-,thick] (2,4) to (-2,6);

\node at (0,0) {$\bullet$};
\node at (2,2) {$\bullet$};
\node at (-2,2) {$\bullet$};
\node at (2,4) {$\bullet$};
\node at (-2,4) {$\bullet$};
\node at (2,6) {$\bullet$};
\node at (-2,6) {$\bullet$};
\node at (0,8) {$\bullet$};
\end{scope}

\begin{scope}[xshift=22.5cm, yshift=15cm]

\draw[-, thick] (0,0) to (-2,2) to (-2,4) to (0,8) to (0,6) to (2,4) to (2,2) to (0,0);
\draw[-, thick] (2,2) to (-2,4);
\draw[-, thick] (-2,2) to (0,6);
\node at (0,0) {$\bullet$};
\node at (2,2) {$\bullet$};
\node at (-2,2) {$\bullet$};
\node at (2,4) {$\bullet$};
\node at (-2,4) {$\bullet$};
\node at (0,6) {$\bullet$};
\node at (0,8) {$\bullet$};
\end{scope}


\begin{scope}[xshift=-22.5cm, yshift=-2cm]

\draw[-,thick] (0,0) to (-2,2) to (-3,6) to (0,8) to (3,6) to (2,2) to (0,0);
\draw[-,thick] (-2,2) to (-2,4) to (3,6) ;
\draw[-,thick] (2,2) to (2,4) to (-3,6) ;
\node at (0,0) {$\bullet$};
\node at (-2,2) {$\bullet$};

\node at (2,2) {$\bullet$};
\node at (2,4) {$\bullet$};
\node at (-2,4) {$\bullet$};
\node at (-3,6) {$\bullet$};
\node at (3,6) {$\bullet$};
\node at (0,8) {$\bullet$};

\end{scope}

\begin{scope}[xshift=-15cm, yshift=-2cm]
\draw[-,thick] (0,0) to (-2,2) to (-2,4) to (0,6) to (2,4) to (2,2) to (0,0);
\draw[-,thick] (-2,2) to (2,4);
\draw[-,thick] (2,2) to (-2,4);
\draw[-,thick] (0,6) to (0,8);
\node at (0,0) {$\bullet$};
\node at (2,2) {$\bullet$};
\node at (-2,2) {$\bullet$};
\node at (2,4) {$\bullet$};
\node at (-2,4) {$\bullet$};
\node at (0,6) {$\bullet$};
\node at (0,8) {$\bullet$};
\end{scope}

\begin{scope}[xshift=-7.5cm, yshift=-2cm]
\draw[-,thick] (0,0) to (-2,2) to (-2,4) to (0,6) to (2,4) to (2,2) to (0,0);
\draw[-,thick] (-2,2) to (2,4);
\draw[-,thick] (2,2) to (-2,4);
\draw[-,thick] (0,0) to (0,-2);
\node at (0,0) {$\bullet$};
\node at (2,2) {$\bullet$};
\node at (-2,2) {$\bullet$};
\node at (2,4) {$\bullet$};
\node at (-2,4) {$\bullet$};
\node at (0,6) {$\bullet$};
\node at (0,-2) {$\bullet$};
\end{scope}

\begin{scope}[xshift=0cm, yshift=-2cm]

\draw[-, thick] (0,0) to (-2,2) to (-2,4) to (0,8) to (0,6) to (2,4) to (2,2) to (0,0);
\draw[-, thick] (2,2) to (-2,4);
\draw[-, thick] (-2,2) to (0,6);
\draw[-, thick] (0,8) to (0,10);
\node at (0,0) {$\bullet$};
\node at (2,2) {$\bullet$};
\node at (-2,2) {$\bullet$};
\node at (2,4) {$\bullet$};
\node at (-2,4) {$\bullet$};
\node at (0,6) {$\bullet$};
\node at (0,8) {$\bullet$};
\node at (0,10) {$\bullet$};
\end{scope}

\begin{scope}[xshift=7.5cm, yshift=-2cm]
\draw[-,thick] (0,0) to (-2,2) to (-2,4) to (0,6) to (2,4) to (2,2) to (0,0);
\draw[-,thick] (-2,2) to (2,4);
\draw[-,thick] (2,2) to (-2,4);
\draw[-,thick] (0,0) to (0,-2);
\draw[-,thick] (0,6) to (0,8);
\node at (0,0) {$\bullet$};
\node at (2,2) {$\bullet$};
\node at (-2,2) {$\bullet$};
\node at (2,4) {$\bullet$};
\node at (-2,4) {$\bullet$};
\node at (0,6) {$\bullet$};
\node at (0,-2) {$\bullet$};
\node at (0,8) {$\bullet$};
\end{scope}

\begin{scope}[xshift=15cm, yshift=-2cm]

\draw[-, thick] (0,0) to (-2,2) to (-2,4) to (0,8) to (0,6) to (2,4) to (2,2) to (0,0);
\draw[-, thick] (2,2) to (-2,4);
\draw[-, thick] (-2,2) to (0,6);
\draw[-, thick] (0,0) to (0,-2);
\node at (0,0) {$\bullet$};
\node at (2,2) {$\bullet$};
\node at (-2,2) {$\bullet$};
\node at (2,4) {$\bullet$};
\node at (-2,4) {$\bullet$};
\node at (0,6) {$\bullet$};
\node at (0,8) {$\bullet$};
\node at (0,-2) {$\bullet$};
\end{scope}

\begin{scope}[xshift=22.5cm, yshift=-2cm]

\draw[-, thick] (0,0) to (-2,2) to (-2,4) to (0,8) to (0,6) to (2,4) to (2,2) to (0,0);
\draw[-, thick] (2,2) to (-2,4);
\draw[-, thick] (-2,2) to (0,6);
\draw[-, thick] (0,0) to (0,-2);
\draw[-, thick] (0,8) to (0,10);
\node at (0,0) {$\bullet$};
\node at (2,2) {$\bullet$};
\node at (-2,2) {$\bullet$};
\node at (2,4) {$\bullet$};
\node at (-2,4) {$\bullet$};
\node at (0,6) {$\bullet$};
\node at (0,8) {$\bullet$};
\node at (0,10) {$\bullet$};
\node at (0,-2) {$\bullet$};
\end{scope}

\end{tikzpicture}
}
\caption{The inclusion posets of $SL(n)$ orbit closures in complexity 1.}
\label{F:28}
\end{figure}

Once again, by the appropriate reduction arguments as we did in the previous cases, 
we may assume that $q_1,q_2,q_3\leq 2$. 
Then we obtain the following five cases: 

\begin{enumerate}
\item $Bl(P_I)=(1,1,2)$ and $Bl(P_J)=(1,1,2)$;
\item $Bl(P_I)=(1,1,2)$ and $Bl(P_J)=(1,2,1)$;
\item $Bl(P_I)=(1,1,2)$ and $Bl(P_J)=(2,1,1)$;
\item $Bl(P_I)=(1,1,3)$ and $Bl(P_J)=(1,2,2)$;
\item $Bl(P_I)=(1,1,3)$ and $Bl(P_J)=(2,1,2)$;
\item $Bl(P_I)=(1,1,3)$ and $Bl(P_J)=(2,2,1)$;
\item $Bl(P_I)=(1,1,4)$ and $Bl(P_J)=(2,2,2)$;
\end{enumerate}
The Hasse diagrams of the corresponding posets of these seven cases 
are given by $P.23$, $P.21$, $P.24$, $P.25$, $P.26$, $P.27$, and $P.28$ in Figure~\ref{F:28}.

\newpage

\section{Proof of Theorem~\ref{T:prelim}}\label{S:Conclusions}

Let $P_I$ and $P_J$ be two standard parabolic subgroups of $SL(n)$ such that 
the diagonal action $SL(n): SL(n)/P_I\times SL(n)/P_J$ has complexity 1. 
Then the block sizes of $P_I$ and $P_J$'s are listed in Table~\ref{T:Ponomareva}.
Let $P$ denote the corresponding inclusion poset of the $SL(n)$-orbit closures. 
The computations that we performed in the previous section show 
that the Hasse diagram of $P$ is one of the 28 non-isomorphic Hasse diagrams 
which are depicted in Figure~\ref{F:28}.
Using this figure, it is easy to verify the following assertions:
\begin{itemize}
\item the posets $P.i$ ($i\in \{1,\dots, 28\}$) have at most 10 elements;
\item $P.21, P.22, P.25, P.27, P.28$ are non-graded posets. 
The poset $P.21$ appears in both of the cases of the 7-th and the 8-th rows of Table~\ref{T:Ponomareva}.
The poset $P.22$ appears only in the case of the 7-th row of Table~\ref{T:Ponomareva}.
The posets $P.25, P.27$, and $P.28$ appear only in the case of the 8-th row of Table~\ref{T:Ponomareva}.
\item The posets $P.1- P.20$ are lattices, and $P.21-P.28$ are non-lattices. 
In particular, all posets of the 8-th row of Table~\ref{T:Ponomareva}
are non-lattices, and two of the five posets of the 7-th row, namely $P.21$ and $P.22$, are non-lattices. 
\end{itemize}
This finishes the proof of Theorem~\ref{T:prelim}.

\begin{Remark}
The {\em height of a finite poset} is the maximum of the lengths of its saturated chains. 
The maximum of the set of the heights of $P.i$'s ($i\in \{1,\dots, 28\}$) is 6.
\end{Remark}

\section{Final Remarks}\label{S:Final}

There is an alternative approach, which is attributed to Bongartz~\cite{Bongartz}, 
for studying the order relations between the closures of the diagonal $G$-orbits in 
a double flag variety \hbox{$G/P_I\times G/P_J$}.
For completeness and for the convenience to the reader, next, 
we will summarize this approach, as described by Magyar, Weyman, and Zelevinsky in~\cite[Example 4.7]{MWZ1};
the exact statement is somewhat difficult to see from Bongartz's original paper.

Let $P_I$ and $P_J$ denote the corresponding standard parabolic subgroups in $G$. 
As before, let us denote by $Bl(P_I)=(p_1,\dots, p_r)$ and $Bl(P_J)=(q_1,\dots, q_s)$ the sizes of the blocks 
of $P_I$ and $P_J$, respectively. For every $G$-orbit in $G/P_I\times G/P_J$, there is a nonnegative 
integer matrix $M=(m_{ij})$ with row sums $p_1,\dots, p_r$, and with column sums $q_1,\dots, q_s$. 
By~\cite[Proposition 4.5]{MWZ1} and by the general results of Bongartz, 
if $M=(m_{ij})$ and $M'=(m'_{ij})$ are two such matrices corresponding to the $G$-orbits $O_1$ and $O_2$ 
in $G/P_I\times G/P_J$, then 
\begin{align}\label{A:notsoobvious}
O_1 \subseteq \overline{O_2} \iff \sum_{k=1}^i \sum_{l=1}^j m_{kl} \geq \sum_{k=1}^i \sum_{l=1}^j m'_{kl}
\ \text{ for all $i$ and $j$.}
\end{align}
This ordering is helpful if the data of two $G$-orbits are provided. 
However, to obtain the full Hasse diagram from (\ref{A:notsoobvious}),  
one needs to generate all possible matrices, and then compare them by using the double-summations as in (\ref{A:notsoobvious}).  
We checked our Figure~\ref{F:28} by following these steps as well.
Our conclusion is that the amount of work that is required for the creation of the Hasse diagrams in our method, 
which uses the minimal double coset representatives, 
and the method of Bongartz, which uses matrices, do not significantly differ from each other. 

\vspace{.25cm}

\textbf{Acknowledgements.}
The first author is partially supported by a grant from
Louisiana Board of Regents. 
The authors are grateful to John Stembridge for his publicly available software codes which were used in the computations of this paper.
The authors thank the referees for their very helpful comments and suggestions.

\bibliography{References}

\begin{thebibliography}{10}

\bibitem{BinghamCanOzan}
Aram Bingham, Mahir~Bilen Can, and Yildiray Ozan.
\newblock A filtration on equivariant {B}orel-{M}oore homology.
\newblock {\em Forum Math. Sigma}, 7:e18, 13, 2019.

\bibitem{Bongartz}
Klaus Bongartz.
\newblock On degenerations and extensions of finite-dimensional modules.
\newblock {\em Adv. Math.}, 121(2):245--287, 1996.

\bibitem{Borel}
Armand Borel.
\newblock {\em Linear algebraic groups}, volume 126 of {\em Graduate Texts in
  Mathematics}.
\newblock Springer-Verlag, New York, second edition, 1991.

\bibitem{Brion89}
Michel Brion.
\newblock Groupe de {P}icard et nombres caract\'{e}ristiques des
  vari\'{e}t\'{e}s sph\'{e}riques.
\newblock {\em Duke Math. J.}, 58(2):397--424, 1989.

\bibitem{Can18}
Mahir~Bilen Can.
\newblock The cross-section of a spherical double cone.
\newblock {\em Adv. in Appl. Math.}, 101:215--231, 2018.

\bibitem{Can19}
Mahir~Bilen Can.
\newblock On the dual canonical monoids, 2019.
\newblock arXiv:1905.08316.

\bibitem{Curtis85}
Charles~W. Curtis.
\newblock On {L}usztig's isomorphism theorem for {H}ecke algebras.
\newblock {\em J. Algebra}, 92(2):348--365, 1985.

\bibitem{DeligneLusztig}
Pierre~R. Deligne and George Lusztig.
\newblock Representations of reductive groups over finite fields.
\newblock {\em Ann. of Math. (2)}, 103(1):103--161, 1976.

\bibitem{DigneMichel}
Fran\c{c}ois Digne and Jean~C.M. Michel.
\newblock Parabolic {D}eligne-{L}usztig varieties.
\newblock {\em Adv. Math.}, 257:136--218, 2014.

\bibitem{HohlwegSkandera}
Christophe Hohlweg and Mark Skandera.
\newblock A note on {B}ruhat order and double coset representatives.
\newblock https://arxiv.org/abs/math/0511611, 2005.

\bibitem{Littelmann}
Peter Littelmann.
\newblock On spherical double cones.
\newblock {\em J. Algebra}, 166(1):142--157, 1994.

\bibitem{MWZ1}
Peter Magyar, Jerzy Weyman, and Andrei Zelevinsky.
\newblock Multiple flag varieties of finite type.
\newblock {\em Adv. Math.}, 141(1):97--118, 1999.

\bibitem{MWZ2}
Peter Magyar, Jerzy Weyman, and Andrei Zelevinsky.
\newblock Symplectic multiple flag varieties of finite type.
\newblock {\em J. Algebra}, 230(1):245--265, 2000.

\bibitem{Panyushev99}
Dmitri~I. Panyushev.
\newblock Complexity and rank of actions in invariant theory.
\newblock {\em J. Math. Sci. (New York)}, 95(1):1925--1985, 1999.
\newblock Algebraic geometry, 8.

\bibitem{Ponomareva}
Elizaveta~V. Ponomareva.
\newblock Classification of double flag varieties of complexity 0 and 1.
\newblock {\em Izv. Ross. Akad. Nauk Ser. Mat.}, 77(5):155--178, 2013.

\bibitem{Stembridge}
John~R. Stembridge.
\newblock Multiplicity-free products and restrictions of {W}eyl characters.
\newblock {\em Represent. Theory}, 7:404--439, 2003.

\bibitem{Therkelsen_Thesis}
Ryan~K. Therkelsen.
\newblock {\em The conjugacy poset of a reductive monoid}.
\newblock ProQuest LLC, Ann Arbor, MI, 2010.
\newblock Thesis (Ph.D.)--North Carolina State University.

\bibitem{Timashev00}
Dmitry~A. Timashev.
\newblock Cartier divisors and geometry of normal {$G$}-varieties.
\newblock {\em Transform. Groups}, 5(2):181--204, 2000.

\bibitem{Timashev}
Dmitry~A. Timashev.
\newblock {\em Homogeneous spaces and equivariant embeddings}, volume 138 of
  {\em Encyclopaedia of Mathematical Sciences}.
\newblock Springer, Heidelberg, 2011.
\newblock Invariant Theory and Algebraic Transformation Groups, 8.

\end{thebibliography}
\bibliographystyle{plain}

\end{document}